\documentstyle[11pt,fleqn]{article}
                      
\title{\bf{ Specht Modules for Weyl Groups}} 

\author{by} 
\date{{\bf Sait ~Hal\i c\i o{\u g}lu}  and  {\bf A O Morris}}
\vspace{1cm}
\renewcommand{\Box}{\rule{2mm}{2mm}}
\renewcommand{\Box}{\rule{2mm}{2mm}}

\begin{document}
\maketitle
\section{Introduction}
Over fields of characteristic zero, there are well known construction 
of the irreducible representations and of irreducible 
modules, called Specht modules for the symmetric groups $S_{n}$ 
which are based on elegant combinatorial 
concepts connected with Young tableaux etc.(see, e.g.[10]).  
James [8] extended these ideas to construct irreducible 
representations and modules over arbitrary field. Al-Aamily,   
Morris and Peel [ 1 ] showed how this 
construction could be extended to deal with the  Weyl groups of type 
$B_{n}$. In [11] the second author described a possible extension of 
James' work for Weyl groups in general, where Young tableaux are 
interpreted 
in terms of root systems. 
We now modify these results and give an alternative generalisation 
of James' work which 
is an improvement and extension of the original 
approach suggested by Morris.

\section{Some General Results On Weyl Groups}
In this section we establish the notation and state some results on Weyl
 groups which are required later . Standard references for this material
 are N Bourbaki [ 3 ] , R W Carter [ 4 ] and J E Humphreys [ 7 ] .

Let  $\Phi$ be a root system in an $l$-dimensional real space $V$ and 
$\pi$  be a simple system in $\Phi$ with corresponding positive system 
$\Phi^{+}$ and negative system $\Phi^{-}$ . 
\[ {\cal W}={\cal W}(\Phi)=<\tau_{\alpha}\mid \alpha \in \Phi >\]
\noindent
be the $Weyl~group$ of $\Phi$, where $\tau_{\alpha}$ is the reflection
corresponding to $\alpha$. Let $l(w)$ denote the $length$ of $w$ 
and the $sign$ of $w$, 
$s(w)$, is defined by 
$s(w)=(~-1~)^{l(w)}~,~w \in {\cal W}$.

We note the following facts which are requried later.

\noindent
{\bf 2.1} There are $\mid {\cal W} \mid$ simple systems (positive systems) 
in $\Phi$ 
given by $w \pi $ $( w \Phi^{+})$, $w \in {\cal W}$. The group ${\cal W}$ 
acts transitively on the set of simple systems.

\noindent
{\bf 2.2} Let $\Gamma$ be the $Dynkin$ $diagram$ ( or $Coxeter~graph$ ) of 
$\Phi$. A Weyl group is $irreducible$ if its Dynkin diagram is connected. 
Irreducible Weyl groups have been classified and correspond to 
root systems of type $A_{l}(l\geq 1)$, $ B_{l} (l\geq 2) $, 
$C_{l} (l\geq 3)$, $ D_{l} (l \geq 4)$, $ ~E_{6}~, ~E_{7}~,~ E_{8} ~$,$~F_{4}~,~G_{2} $  . 
For example $ {\cal W}~(~A_{l}~)~\cong~S_{l+1} $, 
the symmetric group on the set $\{ 1 , 2 , ... , l + 1  \}$ . 
As our aim in this paper is to generalise ideas from the symmetric groups, 
the root system and Dynkin diagram are given in this case. 
The Dynkin diagram is 

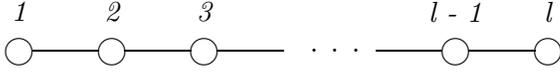
\begin{picture}(350,50)(0,0)   %0,0 is the picture origin
  \put(20,20){\circle{10}}   %20,50 is the centre, 10 is the diameter
  \put(25,20){\line(1,0){25}}  %horizontal vector
   \put(55,20){\circle{10}}
   \put(60,20){\line(1,0){25}} 
   \put(90,20){\circle{10}}
   \put(95,20){\line(1,0){25}}
   \put(130,20){. . .}
   \put(155,20){\line(1,0){25}}
   \put(185,20){\circle{10}}
   \put(190,20){\line(1,0){25}}
   \put(220,20){\circle{10}}
   \put(20,35){\makebox(0,0){\it{1}}}
   \put(55,35){\makebox(0,0){\it{2}}}
   \put(90,35){\makebox(0,0){\it{3}}}
   \put(185,35){\makebox(0,0){\it{l - 1}}}
   \put(220,35){\makebox(0,0){\it{l}}}
 \end{picture}

\noindent
and if $\{~\epsilon_{1}~,~\epsilon_{2}~,~...~,~\epsilon_{l+1}~\}$ is the standard basis for $R^{l+1}$ , then 

\[\pi ~=~\{~\alpha_{1}~=~\epsilon_{1}-\epsilon_{2}~,~\alpha_{2}~=~\epsilon_{2}-\epsilon_{3}~,~...~,~\alpha_{l}~=~\epsilon_{l}- \epsilon_{l+1}~\}\]  
\[\Phi~=~\{~\epsilon_{i}~-~\epsilon_{j}~\mid~1 \leq i ~,~j\leq l+1\}\]
\[\Phi^{+}~=~\{~\epsilon_{i}~-~\epsilon_{j}~\mid~1 \leq i ~<~j\leq l+1\}\]

\noindent
{\bf 2.3} A $subsystem$ $\Psi$ of $\Phi$ is a subset of $\Phi$ which 
is itself a root system in the space which it spans. 
A subsystem $\Psi$ is said to be $additively~closed$ if $\alpha~,~\beta \in \Psi~,~ \alpha~+~\beta \in \Phi~$, then $ \alpha~+~\beta \in \Psi~$ . 
From now on we assume that $\Psi$ is additively closed 
subsystem of $\Phi$. A $Weyl~subgroup$ ${\cal W}(\Psi)$ of ${\cal W}$ corresponding to a 
subsystem $\Psi$ is the subgroup of ${\cal W}$ generated by the 
$\tau_{\alpha} ~,~\alpha \in \Psi $. 

\noindent
{\bf 2.4} The graphs which are Dynkin diagrams of subsystems 
of $\Phi$ may be obtained up to conjugacy by a standard 
algorithm due independently to E B Dynkin, A Borel and J de Siebenthal (see
 e.g. [ 4 ] ).

\noindent
{\bf 2.5}   If $w \in {\cal W}$ and 
$U$ is the subspace of $V$ composed of all vectors fixed by $w$, 
then $w$ is a product of reflections corresponding to roots in 
the orthogonal complement $U^{\perp}$ of $U$ . [ 4 ]

\noindent
{\bf 2.6}  The simple system $J$ of $\Psi$ can always be chosen  such that $J\subset 
\Phi^{+}$. [12]

\noindent
{\bf 2.7} The set $D_{\Psi}=\{ w \in {\cal W} \mid w (j) \in \Phi^{+}~ 
for~all~j \in J~\}$ is a 
$distinguished$ $ set$  $of$ $coset$ $representatives$ of ${\cal W}(\Psi)$ in ${\cal W}$ , 
that is, each element $w \in {\cal W}$ has unique expression of 
the form $d_{\Psi}w_{\Psi}$, where $d_{\Psi} \in D_{\Psi}$ and 
$w_{\Psi} \in {\cal W}(\Psi)$ and furthermore $d_{\Psi}$ is 
the element of minimal length in the coset $d_{\Psi}{\cal W}(\Psi)$. [12]

\noindent
{\bf 2.8} We now extend the above to cover some further reflection subgroups
of Weyl groups. This work is due to Steinberg and we follow the description
 given by Carter [5] where the results highlighted below are proved.

\noindent
Let $\rho$ be a non-trivial symmetry of the Dynkin diagram of $\Phi$. 
Then there is a unique isometry $\tau$ of 
$V$  such that $\tau(r)$ is a positive multiple of $\rho~(~r~)=\bar{r}$ 
for all simple roots $r \in \pi$. The isometry $\tau$ satisfies the conditions :
\[\tau(r)~=~\bar{r}~~~~~~if~ all~ the ~roots~ of ~\Phi ~have ~same~ length\]
\[\tau(r)~=~\{^{{\textstyle \frac{~1}{\sqrt{2}}~\bar{r}~~~if~ r ~is~short}}
_{{\textstyle ~\sqrt{2}~\bar{r}~~if~ r ~is~long}} ~~~ for ~\Phi=B_{2}~ or~ F_{4}\]
\[\tau(r)~=~\{^{{\textstyle \frac{~1}{\sqrt{3}}~\bar{r}~~~if~ r ~is~short}} 
_{{\textstyle ~\sqrt{3}~\bar{r}~~if~ r ~is~long}} ~~~ for ~\Phi=G_{2}\]

Clearly the order of $\tau$ as an isometry of $V$  is equal to the order of $\rho$ as a permutation of $\pi$. Now let

$V^{1}=\{ v \in V \mid \tau(v) = v \}$

\noindent
For each $ v \in V$ let $v^{1}$ denote the projection of $v$ onto the 
subspace $V^{1}$ . Then $v^{1}$ is the average of the vectors in the orbit of $v$ under $\tau$.

Since $\tau(r)$ is a positive multiple of $\rho~(~r~)$ for all $r \in \pi$ we have

$\tau~w_{r}~\tau^{-1}~=~w_{\rho~(~r~)}~,~~~where~ r \in \pi$

Let

${\cal W}^{1}=\{w \in {\cal W} \mid \tau~w~\tau^{-1}~=~w~\}$

\noindent
then ${\cal W}^{1}$ operates faithfully on $V^{1}$. Let $J$  be an orbit of $\pi$ under $\rho$ . Let ${\cal W}(J)$ be the subgroup of ${\cal W}$ generated by the 
elements $\tau_{r}$
 for $r \in J$ . Let $w_{0}^{J}$ be the element of ${\cal W}(J)$
which transforms 
every positive root in $\Phi_{J}$ into a negative root.  
Then $w_{0}^{J} \in {\cal W}^
{1}$ and ${\cal W}^{1}$ is generated by the elements $w_{0}^{J}$ for the different $\rho$-orbits of $\pi$. Since $w_{0}^{J}$ coincides with $w_{r^{1}}$ on 
$ V^{1}$ for each root $r \in J$
the elements $w_{0}^{J}$ are reflections when restricted to $V^{1}$. Then the reflections $w_{r^{1}}$ of $V^{1}$ , for all $r \in \pi$ , generate the group ${\cal W}^{1}$ of isometries of $V^{1}$. Let $\Phi^{1} = \{ r^{1} \in V^{1} \mid r \in \Phi \}$ and $\pi^{1} = \{r^{1} \in V^{1} \mid r \in \pi \}$. The sets $w~(\Phi_{J}^{+})$ form a partition of $\Phi$ as  $w$ runs through 
the elements of ${\cal W}^{1}$ and $J$ runs through the $\rho$-orbits of $\pi$ . The roots $r$ and $s$ are in the same set if and only if $r^{1}$ is a positive multiple of $s^{1}$.

The sets $\Phi^{1}$ and $\pi^{1}$ almost act as a root system 
and simple system for ${\cal W}^{1}$ acting on $V^{1}$, since
\newline
(i) $\Phi^{1}$ spans $V^{1}$ ,
\newline
(ii) Every element of $\Phi^{1}$ is a linear combination of elements of $\pi^{1}$ with coefficients all non-negative or all non-positive ,
\newline
(iii) A basis of $V^{1}$ may be obtained by picking one element of $\pi^{1}$ out of each set of positive multiples ,
\newline
(iv) If $r^{1} \in \Phi^{1}$ then there is an element of ${\cal W}^{1}$ which coincides with $w_{r^{1}}$ on $V^{1}$ ,
\newline
(v) If $r^{1},s^{1} \in \Phi^{1}$ then $w_{r^{1}}(s^{1}) \in \Phi^{1}$.
\newline
The groups ${\cal W}^{1}$ will be referred to as Steinberg subgroups 
from now on.

\section{Specht Modules for Weyl Groups}
Let $\Phi$ be a root system with simple system $\pi$ and Dynkin diagram $\Gamma$
and let $\Psi$ be a subsystem of $\Phi$ with simple system $J \subset \Phi^{+}$
 and  Dynkin diagram $\Delta$.
If  $\Psi= \displaystyle \bigcup_{i=1}^{r}\Psi_{i} $ , 
where $\Psi_{i}$ are the indecomposable components of $\Psi$ , 
then let  $J_{i}$ be a simple system in $\Psi_{i}$ ($i=1,2,...,r$) 
and $J = \displaystyle \bigcup_{i=1}^{r}J_{i}$. Let $\Psi^{\perp}$ be the  largest
subsystem in $\Phi$ orthogonal to $\Psi$ and let $J^{\perp} \subset 
\Phi ^{+}$ the simple system of $\Psi^{\perp}$ .

Let  $\Psi^{'}$ be a subsystem of $\Phi$  which is contained in $\Phi \setminus \Psi$ , with simple system $J^{'} \subset \Phi^{+}$ and  Dynkin diagram $\Delta^{'}$. If  $\Psi^{'}=\displaystyle \bigcup_{i=1}^{s}\Psi_{i}^{'} $ , 
where $\Psi_{i}^{'}$ are the indecomposable components of $\Psi^{'}$  
then let  $J_{i}^{'}$ be a simple system in $\Psi_{i}^{'}$ ($i=1,2,...,s$) 
and $J = \displaystyle \bigcup_{i=1}^{s}J_{i}^{'}$. Let $\Psi^{'^{\perp}}$ be the  largest
subsystem in $\Phi$ orthogonal to $\Psi^{'}$ and let $J^{'^{\perp}} 
\subset \Phi ^{+}$ the simple system of $\Psi^{'^{\perp}}$ .

Let $\bar{J}$ stand for the ordered set $\{J_{1},J_{2},...,
J_{r};J_{1}^{'},J_{2}^{'},...,J_{s}^{'}\}$  , 
where in addition the elements in each 
$J_{i}$ and $J_{i}^{'}$ are ordered. Let 
\[{\cal T}_{J,J^{'}}=\{ w\bar{J} \mid w\in {\cal W} \} \]

Now , we  consider under what conditions the elements in the set 
${\cal T}_{J,J^{'}}$ are distinct. We obtain the following lemma.

\noindent
{\bf Lemma 3.1} {\it  $\mid {\cal T}_{J,J^{'}}\mid=\mid {\cal W} \mid$ if and only if ${\cal W}(J^{\perp}) \cap {\cal W}(J^{'^{\perp}})=<~e~>$ .}

\noindent
{\bf Proof} Suppose that $w_{1}\bar{J}=w_{2}\bar{J}$ , where $w_{1},w_{2} \in {\cal W}$ . 
Then $w_{1} J=w_{2}J$ , $w_{1}J^{'}=w_{2}J^{'}$ and since the elements of each component of $J$ and $J^{'}$ are also ordered it follows that 
$w_{2}^{-1}w_{1}(\alpha)=\alpha$ for all $\alpha \in J$ and $w_{2}^{-1}w_{1}(\alpha)=\alpha$ for all $\alpha \in J^{'}$ , 
that is , by ( 2.5 ) $w_{2}^{-1}w_{1} \in {\cal W}(J^{\perp}) \cap 
{\cal W}(J^{'^{\perp}})$. Thus we have that the elements of 
${\cal T}_{J,J^{'}}$ are distinct if and only 
if ${\cal W}(J^{\perp})\cap {\cal W}(J^{'^{\perp}})=<e>$.

Now we can give our principal definition .

\noindent
{\bf Definition 3.2} Let $\Psi$ and $\Psi^{'}$ be subsystems of $\Phi$ with 
simple systems $J$ and $J^{'}$ respectively such that 
$\Psi^{'} \subseteq \Phi \setminus \Psi$ and $J \subset \Phi^{+}$, 
$J^{'} \subset \Phi^{+}$. The pair $\{J,J^{'}\}$ is called a 
$useful ~system $ in $\Phi$ if ${\cal W}(J) \cap {\cal W}(J^{'})=<e>$ and 
${\cal W}(J^{\perp}) \cap {\cal W}(J^{'^{\perp}})=
<e>$. $~~~~\Box$

\noindent
{\bf Remark 1} If $\{J,J^{'}\}$ is a useful system in $\Phi$ , 
then $\{wJ,wJ^{'}\}$ is also a useful system in $\Phi$, 
for $w \in {\cal W}$. Thus , from now on ${\cal T}_{J,J^{'}}$ will be 
denoted by ${\cal T}_{\Delta}$. 

\noindent
{\bf Remark 2} If $\{J,J^{'}\}$ is a useful system in $\Phi$ 
then $\Psi \cap \Psi^{'} = \emptyset$ and $\Psi^{\perp} \cap 
\Psi^{'^{\perp}} = \emptyset$ However the converse is not true in general. 

\noindent
{\bf Definition 3.3} Let $\{J,J^{'}\}$ be a useful ~system  in $\Phi$ . 
Then the elements of ${\cal T}_{\Delta}$ are called $\Delta-tableaux$,  the  
$J$ and $J^{'}$ are called  the $rows$ and the $columns$ of $\{J,J^{'}\}$ 
respectively.$~~~~\Box$

We see that this is a natural extension of the concept of a Young tableaux in the representation theory of symmetric 
groups . In the case $A_{l}$ , the subsystems are given up 
to conjugacy by root systems of type $A_{l_{1}}+...+~A_{l_{s}}$, where the $l_{i} \geq 0$ and $\sum (l_{i}+1~)~=~l+1$   Thus the corresponding Weyl subgroup is ${\cal W}~(A_{l_{1}}~)\times~...~\times~{\cal W}~(~A_{l_{s}}~)$ which is isomorphic to the Young 
subgroup $S_{l_{1}+1}~\times~...~\times~S_{l_{s}+1} $ of $S_{l+1}$ and there is a Weyl subgroup corresponding to each partition of $l~+~1$ . If $\Psi = A_{l_{1}}~+~...~+~A_{l_{s}}$ is 
a subsystem of $A_{l}$ , then $\lambda = (l_{1}+1~,~...~,~l_{s}+1)$ is a 
partition of $l+1$ . If we put $k_{i}=l_{1}+...+l_{i}+i~(i=1,2,...,s)$ and $k_{0}=0$ then $\{\epsilon_{k_{i-1}+1}-
\epsilon_{k_{i-1}+2},...,\epsilon_{k_{i}-1}-\epsilon_{k_{i}}\}$ is a 
simple system for $A_{l_{i}}$ . 
This may be represented by the $\lambda$-tableau

$\begin{array}{cccccccccc}
1&2&3&.&.&.&.&.&.&k_{1}\\
k_{1}+1&k_{1}+2&k_{1}+3&.&.&.&.&.&k_{2}&\\
k_{2}+1&k_{2}+2&k_{2}+3&.&.&.&.&k_{3}&&\\
.&.&.&.&.&.&&&&\\
k_{s-1}+1&k_{s-1}+2&k_{s-1}+3&.&.&l+1&&&&
\end{array}$

\noindent
the other n! $\lambda-tableau$ being obtained by allowing the 
elements of $S_{n}$ to act on this $\lambda-tableau$.  

The orthogonal subsystem $\Psi^{\perp}$ is the root system determined by the elements in rows of length one in the $\lambda-tableau$ . The ordering of the rows is important , for example the $3^{2}-tableaux$ $^{123}_{456}$ and $^{456}_{123}$ are regarded as distinct tableaux.

\noindent
{\bf Definition 3.4}  Two $ \Delta$-tableaux $\bar{J}$ and $\bar{K}$  
are $row-equivalent$ ,  written $\bar{J}~\sim~\bar{K}$ ,  if there exists  $w \in {\cal W}(J)$ such that $\bar{K}~=~w~\bar{J}$.
The equivalence class which contains the $\Delta$-tableau $\bar{J}$ is 
$\{\bar{J}\}$ and is called a $\Delta-tabloid.~~~~~\Box$

Let $\tau_{\Delta}$ be the set of all $\Delta$-tabloids . 
It is clear that the number of distinct elements in 
$\tau_{\Delta}$ is  $[{\cal W}~:{\cal W}(J)]$  and by ( 2.7 )  we have
\[\tau_{\Delta}=\{ \{~d\bar{J}~\} \mid d\in D_{\Psi} \} \]
\noindent 
We note that if $\bar{J}=\{~J~;~J^{'}~\}$ then $dJ \subset \Phi^{+}$ but $dJ^{'}$ need not be a subset of $\Phi^{+}$ .

We now give an example to illustrate the construction of a $\Delta$-tabloid. 
In this example and later examples we use the following notation. 
If $\pi = \{\alpha_{1},\alpha_{2},...,\alpha_{n}\}$ is a simple system in $\Phi$ and $\alpha \in \Phi$, then $\alpha= {\displaystyle \sum_{i=1}^{n}a_{i}\alpha_{i}}$, where $a_{i} \in {\bf Z}$. From now on $\alpha$ is 
denoted by $a_{1}a_{2}...a_{n}$ and $\tau_{\alpha_{1}},
\tau_{\alpha_{2}},...,\tau_{\alpha_{n}}$ are denoted by 
$\tau_{1},\tau_{2},...,\tau_{n}$ respectively. 

\noindent
{\bf Example 3.5} Let $\Phi=\bf{A_{3}}$ with simple system  
$\pi=\{\alpha_{1}=\epsilon_{1}-\epsilon_{2},\alpha_{2}=\epsilon_{2}
-\epsilon_{3},
\alpha_{3}=\epsilon_{3}-\epsilon_{4}\}$.  
 Let $\Psi={\bf 2A_{1}}$ be
 the subsystem of $\Phi$ with simple system 
 $J = \{100,001\}$. Then the Dynkin diagram for $\Psi$ is 
\begin{picture}(145,50)(0,0)   %0,0 is the picture origin
  \put(35,25){\circle{10}}   %20,50 is the centre, 10 is the diameter
 \put(40,25){\line(1,0){25}}  %horizontal vector
 \put(70,25){\makebox(0,0){$\bigotimes$}}
  \put(75,25){\line(1,0){25}} 
  \put(105,25){\circle{10}}
  \put(35,5){\makebox(0,0){\bf{1}}}
  \put(70,5){\makebox(0,0){\bf{2}}}
  \put(105,5){\makebox(0,0){\bf{3}}}
  \end{picture}

\noindent
Let $\Psi^{'}={\bf 2A_{1}}$ be the subsystem of $\Phi$  which is contained in 
$\Phi \setminus \Psi$, with simple system  $J^{'}=\{110,011\}$. 
Since ${\cal W}(J) \cap {\cal W}(J^{'})=<e>$ and   
${\cal W}(J^{\perp}) \cap {\cal W}(J^{'^{\perp}})=<e>$, 
then $\{J,J^{'}\}$ is a useful system in $\Phi$ . 
Then $\tau_{\Delta}$ contains the $\Delta-tabloids$ (we also give 
these in the traditional notation as in [ 9 ] );

\begin {tabular}{llll}
$\{ \bar{J} \}$&=&$\{100,001;110,011\}~$&$~~~~~~~~~(^{\overline {12}}_{{\overline 34}})$
\\
$\{ \tau_{2}\bar{J} \}$&=&$\{110,011;100,001 \}~$&$~~~~~~~~~(^{\overline {13}}_
{\overline{ 24}})$\\
$\{ \tau_{1}\tau_{2}\bar{J }\}$&=&$\{010,111;-100,001 \}~$&$~~~~~~~~~(^
{\overline{ 23}}_{\overline {14}})$\\
$\{ \tau_{3}\tau_{2}\bar{J} \}$&=&$\{111,010;100,-001 \}~$&$~~~~~~~~~(^
{\overline{ 14}}_{\overline{ 23}})$\\
$\{ \tau_{1}\tau_{3}\tau_{2}\bar{J} \}$&=&$\{011,110;-100,-001\}~$&
$~~~~~~~~~(^{\overline{ 24}}_{\overline{ 13}})$\\
$\{ \tau_{2}\tau_{1}\tau_{3}\tau_{2}\bar{J} \}$&=&$\{001,100;-110,-011\}$
&$~~~~~~~~~(^{\overline {34}}_{\overline {12}})~~~~~~~~~\Box$
\end{tabular}

\vspace{0.1cm}
The group ${\cal W}$ acts on $\tau_{\Delta}$ according to 
\[\sigma~\{\overline{wJ}\} =\{\overline{\sigma w J}\}~~~~~~for~all~\sigma \in {\cal W} .\]
\noindent
This action is well defined , for if $\{\overline{w_{1}J}\}=\{\overline{w_{2}J}\}$ , then there exists $\rho \in {\cal W}(w_{1}J)$ such that $\overline{\rho w_{1} J} =\overline{ w_{2}J}$ . 
Hence since $\sigma \rho \sigma^{-1} \in {\cal W}(\sigma w_{1}J)$ and $\overline{\sigma w_{2} J }=\overline{\sigma \rho w_{1} J} = (\sigma \rho \sigma^{-1} )(\overline{\sigma w_{1} J })$ , we have $\{\overline{\sigma w_{1}J}\}=\{\overline{\sigma w_{2}J}\}$ . 

Now if K is arbitrary field , let $M^{\Delta}$ be the 
K-space whose basis elements are the $\Delta$-tabloids. Extend the action of ${\cal W}$ on $\tau_{\Delta}$ linearly on $M^{\Delta}$ , then $M^{\Delta}$ becomes  a $K{\cal W}$-module. Then we have the following lemma.

\noindent
{\bf Lemma 3.6} {\it The $K{\cal W}$-module $M^{\Delta}$  is a cyclic 
$K{\cal W}$-module generated by any one tabloid and 
$dim_{K}M^{\Delta} = [{\cal W}:{\cal W}(J)]~~~~~~\Box$}

Now we proceed to consider the possibility of constructing a $K{\cal W}$-module which corresponds to the Specht module in the case of symmetric groups. In order to do this  we need to define a $\Delta$-polytabloid .

\noindent
{\bf Definition 3.7} Let $\{J,J^{'}\}$ be a useful system in $\Phi$ . Let
\begin {eqnarray*}
\kappa_{J^{'}}~=~\sum_{\sigma~\in~{\cal W}(J^{'})}~s~(~\sigma~)~\sigma~~~~~~and~~~e_{J,J^{'}}~=~\kappa_{J^{'}}~\{~\bar{J}~\}
\end{eqnarray*}
\noindent
where $s$ is the sign function defined in Section 2 . 
Then  $e_{J,J^{'}}$ is called the generalized $\Delta-polytabloid 
$ associated with $J$. $~~~~~~\Box$

If w $\in {\cal W}(\Phi)$ , then 
\begin {eqnarray*}
w~\kappa_{J^{'}}&=&\sum_{\sigma~\in~{\cal W}(J^{'})}~s~(~\sigma~)~w~\sigma\\
&=&\sum_{\sigma~\in~{\cal W}(J^{'})}~s~(~\sigma~)~(w~\sigma~w^{-1})~w\\
&=&~\{~\sum_{\sigma~\in~{\cal W}(wJ^{'})}~s~(~\sigma~)~\sigma~\}~w
\end{eqnarray*}
\noindent
Hence, for all $w \in {\cal W}(\Phi)$, we have

\begin {tabular}{llll}
$w~e_{J,J^{'}}$&=&$\kappa_{wJ^{'}}~\{~\overline{wJ}~\}~=~e_{wJ,wJ^{'}}$~~~~~~~~~~ ~~(~3.1~)&
\end{tabular}

 Let $S^{J,J^{'}}$ be the subspace of $M^{\Delta}$ generated by $e_{wJ,wJ^{'}}$ where $w \in {\cal W}$. Then by (3.1) $S^{J,J^{'}}$ is a 
$K{\cal W}$-submodule of $M^{\Delta}$, 
which is called a $generalized$ $Specht$ $module$. Then we have the following theorem .

\noindent
{\bf Theorem 3.8 } {\it The $K{\cal W}$-module $S^{J,J^{'}}$ is a cyclic submodule generated by any  $\Delta$-polytabloid.}$~\Box$

The following proposition notes some isomorphisms between Specht modules.

\noindent
{\bf Proposition 3.9} {\it Let $\{J,J^{'}\}$ be a useful system in $\Phi$ . Then we have the following isomorphisms:

\begin {tabular}{llll}
$(i)~If~w\in {\cal W}~,$&then&$S^{J,J^{'}}\cong S^{wJ,wJ^{'}}$&\\
$(ii)~If~w\in {\cal W}(J)~,$&then&$S^{J,J^{'}}\cong S^{J,wJ^{'}}$&\\
$(iii)~If~w\in {\cal W}(J^{'})~,$&then&$S^{J,J^{'}}\cong S^{wJ,J^{'}}$&
\end{tabular}}

\noindent
{\bf Proof} (i) If $w\in {\cal W}$ , then define $T_{w}~:~S^{J,J^{'}}~\rightarrow~S^{wJ,wJ^{'}}$ by 
\[T_{w}(~e_{J,J^{'}}~)=w~(~e_{J,J^{'}}~)\]
\noindent
Then by ( 3.1 ) we have $w~(~e_{J,J^{'}}~)=
e_{wJ,wJ^{'}}\in S^{wJ,wJ^{'}}$. Clearly $T_{w}$ is an isomorphism. 
The statements (ii) , (iii)  follow easily from (i).$~~~~\Box$

Proposition 3.9 says  that a generalized Specht module is dependent only 
on the Dynkin diagram $\Delta$ and $\Delta^{'}$ of $J$ and $J^{'}$ 
respectively, thus, from now on it will be denoted by $S^{\Delta,\Delta^{'}}$.

A Specht module is spanned by the $e_{wJ,wJ^{'}}$ for all $ w \in {\cal W} $ ; the next lemma shows that we need only consider a certain subset of ${\cal W}$ .

\noindent
{\bf Lemma 3.10} {\it Let $\{J,J^{'}\}$ be a useful system in $\Phi$ . 
Then $S^{\Delta,\Delta^{'}}$ is generated by $e_{dJ,dJ^{'}}$ ,  
where $d \in D_{\Psi^{'}}$.}

\noindent
{\bf Proof} If $w \in {\cal W}$, then by (2.7), $w = d \rho$, 
where $d \in D_{\Psi^{'}}$ and $\rho \in {\cal W}(J^{'})$ and since 
\begin {eqnarray*}
\rho ~e_{J,J^{'}}&=&\sum_{\sigma~\in~{\cal W}(J^{'})}~s~(~\sigma~)~\rho~\sigma~\{ \bar{J} \} = s~(\rho) ~e_{J,J^{'}} 
\end{eqnarray*}
\noindent
 we have
\[e_{wJ,wJ^{'}} ~= ~w ~e_{J,J^{'}}~=~d~\rho~e_{J,J^{'}}~=~s(\rho)~d~e_{J,J^{'}}~=~s(\rho)~e_{dJ,dJ^{'}}~~~~~~\Box\]

\noindent
{\bf Lemma 3.11} {\it Let $\{J,J^{'}\}$ be a useful system in $\Phi$ and 
let $d \in D_{\Psi}$. If $\{\overline{dJ} \}$ appears in $e_{J,J^{'}}$ then it appears only once.}

\noindent
{\bf Proof} If $\sigma~,~\sigma^{'} \in {\cal W}(J^{'})$ and suppose that $\sigma = dw~,~
\sigma^{'}=dw^{'}$ where $w,w^{'} \in {\cal W}(J) $ . Then $d=\sigma w^{-1}=
\sigma^{'}w^{'^{-1}}$ and $\sigma^{'^{-1}}\sigma=w^{'^{-1}}w \in {\cal W}(J) 
\cap {\cal W}(J^{'})=~<~e~>$ . Hence we have $w=w^{'}$ and $\sigma=\sigma^{'}$ . 
Then $\{\overline{dJ} \}$ appears in $e_{J,J^{'}}$ only once.$~~~~\Box$

\noindent
{\bf Corollary 3.12} {\it  If $\{J,J^{'}\}$ be a useful system in $\Phi$ , then $e_{J,J^{'}} \neq 0$ .}

\noindent
{\bf Proof} By Lemma 3.11 if $\{\overline{\sigma J} \}$ appears in $e_{J,J^{'}}$ 
then  all the $\{\overline{\sigma J }\}$ are different, where $\sigma \in {\cal W}(J^{'})$.  
But $\{~\{\overline{\sigma J }\} \mid \sigma \in {\cal W}(J^{'})~\}$ is a linearly 
independent subset of $\{~\{\overline{d~J }\} \mid d \in D_{\Psi}~\}$.   
If $e_{J,J^{'}} =0$ then $s(\sigma)=(-1)^{l(\sigma)}=0$ for all $\sigma 
\in {\cal W}(J^{'})$ . This is a contradiction and so  $e_{J,J^{'}} \neq 0.$ $~~~\Box$ 

\noindent
{\bf Remark 3} In [ 11 ] the second author  has defined a partial dual of a subsystem $\Psi$ simply as
a subsystem $\Psi^{'}$ of $\Phi$ which is contained in $\Phi\setminus \Psi$ . 
However the following lemma shows that the extra condition 
${\cal W}(J) \cap {\cal W}(J^{'})=<e>$  in our definition of a 
useful system is  also necessary. 
Unfortunately this condition which is a group theoretical one is 
not easily checked and it would be useful if it could be 
replaced by a criteria in terms of the root system only.

\noindent
{\bf Lemma 3.13} {\it If there exists $w \in {\cal W}(J) \cap 
{\cal W}(J^{'})$ such that $w$ has order 2 , and s($w$) = -1 then $e_{J,J^{'}}~=~0$.}

\noindent
{\bf Proof} If $w\in {\cal W}(J)\cap~{\cal W}(J^{'})$ and $w$ has order 2 then
\[(~e-w~)~\{\bar{J}\} = \{\bar{J}\} - \{\bar{J}\} = 0 \]
\noindent
and also $\{~e~,~w~\}$ is a subgroup of ${\cal W}(J^{'})$ . 
Thus we can select signed coset representatives $\sigma_{1}$~,$~\sigma_{2}$~,
~$\sigma_{3}$~,~...~,~$\sigma_{s}$ for $\{~e~,~w~\}$ in ${\cal W}(J^{'})$ such that
\begin {eqnarray*}
e_{J,J^{'}}&=&\sum_{\sigma~\in~{\cal W}(J^{'})}~s~(~\sigma~)~\sigma~\{\bar{J}\}\\
&=&(~\sum_{i~=~1}^{s}~\sigma_{i})~(~e~-~w~)~\{\bar{J}\}\\
&=&0.~~~~~~~~~~~~~~~~~~~~\Box
\end{eqnarray*}

\noindent
{\bf Example 3.14} Let $\Phi={\bf B_{3}}$ and 
$\pi=\{\alpha_{1}=\epsilon_{1}-\epsilon_{2},\alpha_{2}
=\epsilon_{2}-\epsilon_{3},\alpha_{3}=\epsilon_{3} \}$. Let $\Psi=
{\bf 3A_{1}}$ be the subsystem of $\Phi$ with simple system 
$J=\{\alpha_{1}=\epsilon_{1}-\epsilon_{2},\tilde{\alpha}=\epsilon_{1}+\epsilon_{2},\alpha_{3}=\epsilon_{3} \}$ and 
let $\Psi^{'}={\bf 3A_{1}}$ be the subsystem of $\Phi$ with $J^{'}=\{\alpha_{2}=\epsilon_{2}-\epsilon_{3},\alpha_{1}+\alpha_{2}+\alpha_{3}=\epsilon_{1},\alpha_{2}+2~\alpha_{3}=\epsilon_{2}+\epsilon_{3}\}$. Then $\Psi \cap \Psi^{'}=\emptyset$. But
\newline
${\cal W}(J)~~$=$~\{e,\tau_{1},\tau_{3},\tau_{1}\tau_{3},
\tau_{2}\tau_{3}\tau_{1}\tau_{2}\tau_{3}\tau_{1}\tau_{2},\tau_{2}\tau_{3}\tau_{1}\tau_{2}\tau_{3}\tau_{1}\tau_{2}\tau_{1},
\tau_{3}\tau_{2}\tau_{3}\tau_{1}\tau_{2}\tau_{3}\tau_{1}\tau_{2},\tau_{3}\tau_{2}\tau_{3}\tau_{1}\tau_{2}\tau_{3}\tau_{1}\tau_{2}\tau_{1}\}$ 
\newline
${\cal W}(J^{'})~$=$~\{e,\tau_{2},
\tau_{1}\tau_{2}\tau_{3}\tau_{2}\tau_{1},\tau_{3}\tau_{2}\tau_{3},\tau_{3}\tau_{2}\tau_{3}\tau_{1}\tau_{2}\tau_{3}\tau_{1}\tau_{2}\tau_{1},\tau_{3}\tau_{2}\tau_{3}\tau_{2},
\tau_{3}\tau_{2}\tau_{3}\tau_{1}\tau_{2}\tau_{3}\tau_{2}\tau_{1},\tau_{1}\tau_{2}\tau_{3}\tau_{1}\tau_{2}\tau_{1}\}$ 

\noindent 
It follows that $w=~\tau_{3}\tau_{2}\tau_{3}\tau_{1}\tau_{2}\tau_{3}\tau_{1}\tau_{2}
\tau_{1}~ \in {\cal W}(J) \cap {\cal W}(J^{'}) $ and $e_{J,J^{'}}=0$ .$~~~~\Box$

\noindent
{\bf Lemma 3.15} {\it  Let $\{J,J_{1}^{'}\}$ and $\{J,J_{2}^{'}\}$ be useful systems in $\Phi$ . If $\Psi_{1}^{'} \subseteq \Psi_{2}^{'}$ , then $S^{J,J_{2}^{'}}$ is a $K{\cal W}$-submodule of $S^{J,J_{1}^{'}}$ , where $J_{1}^{'}$ and $J_{2}^{'}$  are simple systems for $\Psi_{1}^{'}$ and $\Psi_{2}^{'}$ respectively.}

\noindent
{\bf Proof} Since $\Psi_{1}^{'} \subseteq \Psi_{2}^{'}$ , 
${\cal W}(J_{1}^{'})$ is a subgroup of ${\cal W}(J_{2}^{'})$ and 
we can select coset representatives $a_{1},a_{2},...,a_{n}$ for 
${\cal W}(J_{1}^{'})$ in ${\cal W}(J_{2}^{'})$ such that  
${\cal W}(J_{2}^{'})= \displaystyle \bigcup_{i=1}^{n} a_{i}
{\cal W}(J_{1}^{'})$ . Then we have
\begin {eqnarray*}
e_{J,J_{2}^{'}}&=&\sum_{\sigma~\in~{\cal W}(J_{2}^{'})}~s~(~\sigma~)~\sigma~\{\bar{J}\}\\
&=&~\sum_{i~=~1}^{n}~a_{i}~s(a_{i})\sum_{\sigma~\in~{\cal W}(J_{1}^{'})}~s~(~\sigma~)~\sigma~\{\bar{J}\}\\
&=&~(~\sum_{i~=~1}^{n}~a_{i}~s(a_{i})~)~e_{J,J_{1}^{'}}~~~~~~~~~~~~~~\Box
\end{eqnarray*}

Now we  consider under what conditions  $S^{\Delta,\Delta^{'}}$ is irreducible .

\noindent
{\bf Lemma 3.16} {\it Let $\{J,J^{'}\}$ be a useful system in $\Phi$ and  let $ d \in D_{\Psi}$ . Then the following conditions are equivalent:
\newline
(i) $\{~\overline{dJ}~\}$ appears with non-zero coefficient in $e_{J,J^{'}}$
\newline
(ii) There exists $\sigma \in {\cal W}(J^{'})$ such that $\sigma \{~\bar{J}~\}~=~\{\overline{dJ}~\}$
\newline
(iii) There exists $\rho \in {\cal W}(J)$  and $\sigma \in {\cal W}(J^{'})$ such that $d~=~\sigma~\rho$}

\noindent
{\bf Proof}  The equivalence of (i) and (ii) follows directly from the formula 
\[e_{J,J^{'}}=\sum_{\sigma~\in~{\cal W}(J^{'})}~s~(~\sigma~)~\sigma~\{\bar{J}\}\]
\noindent
(ii) $\Rightarrow$ (iii) Suppose that there exists $\sigma \in 
{\cal W}(J^{'})$ such that $\sigma \{~\bar{J}~\}~=~\{\overline{dJ}~\}$. 
Then we have $\sigma^{-1} ~d\{~\bar{J}~\}~=~\{\bar{J}~\}$. 
By the definition of 
equivalence  there exists $\rho \in {\cal W}(J)$ such that 
$\sigma^{-1} ~d~\bar{J}~=~\rho~\bar{J}~$. Then  $\rho^{-1}\sigma^{-1}d 
\in {\cal W}(J^{\perp}) \cap {\cal W}(J^{'^{\perp}})$. 
Since $\{J,J^{'}\}$ is a useful system in $\Phi$ then $d~=~\sigma~\rho$ , where $\sigma \in {\cal W}(J^{'})$ and $\rho \in {\cal W}(J)$.
\newline
(iii) $\Rightarrow$ (ii) Let $d~=~\sigma~\rho$, where 
$\sigma \in {\cal W}(J^{'})$ and $\rho \in {\cal W}(J)$. 
Since $\rho \in {\cal W}(J)$,  $~\rho~\bar{J}~=~\bar{J}$  
then  $\{\overline{dJ}~\}~=~\{\overline{\sigma \rho J}~\}~=~
\{\overline{\sigma J}~\}.~~~~~\Box $

\noindent
{\bf Lemma 3.17} {\it Let $\{J,J^{'}\}$ be a useful system in $\Phi$ and  let $ d \in D_{\Psi}$ . If $\{~\overline{dJ}~\}$ appears in
 $e_{J,J^{'}}$ then $d~\Psi \cap \Psi^{'}~=~\emptyset $.}

\noindent
{\bf Proof}
 Let $\alpha \in d~\Psi$. If $\{~\overline{dJ}~\}$ appears 
in $e_{J,J^{'}}$ then by Lemma 3.16 $d~=~\sigma~\rho$, where $\sigma \in 
{\cal W}(J^{'})$ and $\rho \in {\cal W}(J)$. Then 
$ \alpha \in \sigma \rho \Psi$. Since $\rho \in {\cal W}(J)$, 
then $ \alpha \in \sigma \Psi$ and $\sigma^{-1}\alpha \in \Psi$. But 
$\Psi \cap \Psi^{'} = \emptyset$, then $\sigma^{-1}\alpha \not \in \Psi^{'}$. 
Since $\sigma \in {\cal W}(J^{'})$, $\sigma \Psi^{'}~=~\Psi^{'}$ then $\alpha \not \in \Psi^{'}$. $~~~~\Box$

\noindent
{\bf Lemma 3.18} {\it Let $\{J,J^{'}\}$  be  a useful system in $\Phi$ and 
let $ d \in D_{\Psi}$. Let \newline
$d~\Psi \cap \Psi^{'} \not = \emptyset$. 
Then $\kappa_{J^{'}}\{~\overline{dJ}~\} = 0 $.}

\noindent
{\bf Proof} If $d~\Psi \cap \Psi^{'} \not = \emptyset$ then let $\alpha \in d~\Psi \cap \Psi^{'} $ and so $\tau_{\alpha} \in
 {\cal W}(dJ) \cap {\cal W}(J^{'})$. Thus 
\[(~e-\tau_{\alpha}~)~\{\overline{dJ}\} =
 \{\overline{dJ}\} - \{\overline{dJ}\} = 0 \]
\noindent
Since $\{~e~,~\tau_{\alpha}~\}$ is a 
subgroup of ${\cal W}(J^{'})$ then we can select signed coset 
representatives $\sigma_{1}$, $~\sigma_{2}$, ~$\sigma_{3}$, ~...~, ~$\sigma_{s}$ for $\{~e~,~\tau_{\alpha}~\}$ in ${\cal W}(J^{'})$ such that

\[\kappa_{J^{'}}\{~\overline{dJ}~\} =(~\sum_{i~=~1}^{s}~\sigma_{i})~(~e~-~\tau_{\alpha}~)~\{\overline{d~J}\}~=~0 ~~. ~~~~~~~~\Box \]

The converse of Lemma 3.17 is not true in general, which leads to the following definition .

\noindent
{\bf Definition 3.19}  A useful system $\{J,J^{'}\}$ in $\Phi$ is 
called a $good~system$  if 
\newline
$d~\Psi \cap \Psi^{'} = \emptyset$ for $d \in D_{\Psi}$ then $\{~\overline{dJ}~\}$ 
appears with non-zero coefficient in $e_{J,J^{'}}$.$~\Box$

\noindent
{\bf Lemma 3.20} {\it  Let $\{J,J^{'}\}$ be a good system in $\Phi$ and  let $d \in D_{\Psi}$ .
\newline
(i) If $\{\overline{dJ}\}$ does not appear in $e_{J,J^{'}}$  then $\kappa_{J^{'}}\{~\overline{dJ}~\} = 0 $.
\newline
(ii) If $\{\overline{dJ}\}$  appears in $e_{J,J^{'}}$  then there exists $ \sigma \in {\cal W}(J^{'})$ such that}
\[\kappa_{J^{'}}\{~\overline{dJ}~\}~=~s~(~\sigma~)~e_{J,J^{'}} \]

\noindent
{\bf Proof}  (i) If $\{\overline{dJ}\}$ does not appear in $e_{J,J^{'}}$  
then by definition of a good system  
we have $d~\Psi \cap \Psi^{'}  \not = \emptyset$ and hence by Lemma 3.17 we have $\kappa_{J^{'}}\{~\overline{dJ}~\} = 0 $.
\newline
(ii) Since $\{\overline{dJ}\}$  appears in $e_{J,J^{'}}$ it follows by 
Lemma 3.16 that there exists $\sigma \in {\cal W}(J^{'})$ such that $\sigma \{\bar{J}\}~=~\{\overline{dJ}\}$. Then we have

\[\kappa_{J^{'}}\{~\overline{dJ}~\} =(~\sum_{\rho \in {\cal W}(J^{'})}~s(\rho)~\rho)~\{\overline{\sigma~J}\}~=~s(\sigma)e_{J,J^{'}} ~~.
~~~~~~~\Box \]

\noindent
{\bf Corollary 3.21} {\it  Let $\{J,J^{'}\}$ be a 
good system in $\Phi$.  If $m \in M^{\Delta}$ then $\kappa_{J^{'}} 
m$ is a
 multiple of $e_{J,J^{'}}$.}

\noindent
{\bf Proof} 
 If $m \in M^{\Delta}$ then we have 
\begin {eqnarray*}
m&=& \sum_{d ~\in ~D_{\Psi}} ~\alpha_{d}
 ~\{ ~\overline{dJ}~\}~, where ~\alpha_{d} \in K \\
\kappa_{J^{'}}~m&=& \sum_{d ~\in ~D_{\Psi}}
 ~\alpha_{d}~\kappa_{J^{'}}~\{ ~\overline{dJ}~\}
\end{eqnarray*}

\noindent
But by Lemma 3.20  $\kappa_{J^{'}}~m~=~
\lambda~e_{J,J^{'}}$ , where $\lambda \in K$.
 $~~~~~~~~~~\Box$

We now define a bilinear form $ <~,~>$ on 
$M^{\Delta}$ by setting
\[<~\{\bar{J_{1}}\}~,~\{\bar{J_{2}}\}~>~=
\{^{~~1~~~~if~\{\bar{J_{1}}\}~=~\{\bar{J_{2}}\}}_
{~~0~~~~otherwise} \]

\noindent
This is a symmetric, non-singular, ${\cal W}$-invariant bilinear form
 on $M^{\Delta}$ .

Now we can prove James' submodule theorem in this general setting.

\noindent
{\bf Theorem 3.22} {\it  Let $\{J,J^{'}\}$ be a good 
system in $\Phi$.  Let U be submodule of $M^{\Delta}$ .
 Then either $S^{\Delta , \Delta^{'}} \subseteq U$ or $U 
\subseteq S^{\Delta ,
 \Delta^{'^\perp}}$ , where $S^{\Delta,\Delta^{'^\perp}}$
 is complement of $S^{\Delta,\Delta^{'}}$ in  $M^{\Delta} .$}

\noindent
{\bf Proof} If $u \in U$ then
\begin {eqnarray*}
<~u~,~e_{J,J^{'}}~>~&=& <~u~,\sum_{\sigma ~\in 
~{\cal W}(J^{'})} ~s~(~\sigma~)~\sigma 
~\{ ~\bar{J}~\}~> \\
&=&~\sum_{\sigma ~\in ~{\cal W}(J^{'})}
 ~<~s~(~\sigma~)~\sigma 
^{-1}~u~,~\{ ~\bar{J}~\}~> \\
&=&<~\kappa_{J^{'}}~u~,~\{ ~\bar{J}~\}~>
\end{eqnarray*}
\noindent
But by Corollary 3.21 $\kappa_{J^{'}}~u~=~
\lambda~e_{J,J^{'}}$ , for some $\lambda \in K$. 
 If $\lambda \not = 0$ for some $u \in U$ , then
 $e_{J,J^{'}} \in U$ , 
that is , $S^{\Delta,\Delta^{'}} \subseteq U$. 
However , if $\lambda  = 0$ for all $u \in U$ , then
 $<~u~,~e_{J,J^{'}}~>~=0$ , that is , $U \subseteq 
S^{\Delta,\Delta^{'^\perp}}$.
$~~~~~~~\Box$

We can now prove our principal result.

\noindent
{\bf Theorem 3.23} {\it Let $\{J,J^{'}\}$ be a
 good system in $\Phi$ . The $K{\cal W}$-module 
\newline
$D^{\Delta,\Delta^{'}}~=~S^{\Delta,\Delta^{'}}~
/~S^{\Delta,\Delta^{'}}\cap S^{\Delta,\Delta^{'^\perp}}$ 
is zero or irreducible.}

\noindent
{\bf Proof} 
 If U is a submodule of $S^{\Delta,\Delta^{'}}$ 
then U is a 
submodule of $M^{\Delta}$ and by Theorem 3.22 
either $S^{\Delta,\Delta^{'}} 
\subseteq U$ in which case $U~=~
S^{\Delta,\Delta^{'}} $ or $U \subseteq 
S^{\Delta,\Delta^{'^\perp}}$ and $U
 \subseteq ~S^{\Delta,\Delta^{'}}~\cap ~
S^{\Delta,\Delta^{'^\perp}}$ , which 
completes the proof. $~~~~~\Box$

In the case of K = {\bf Q} or any field of
 characteristic zero $<~,~>$ is an inner product and
 $D^{\Delta,\Delta^{'}}~=~S^{\Delta,\Delta^{'}}$.
 Thus if for a 
subsystem $\Psi$ of $\Phi$
 a good system $\{ J,J^{'}\}$ can be found , 
then we have a construction for irreducible 
$K{\cal W}$-modules . Hence it is essential to 
show for each subsystem 
that a good system exists which satisfies
 Definition 3.19.

In the following example, we show how a good system 
may be constructed in most cases for the Weyl group of 
type ${\bf D_{4}}$.
 In a future publication, we shall present an algorithm 
for constructing a good system for certain subsystems; 
indeed this
 algorithm will give a good system with additional
 properties which will lead to the construction
 of a $K$-basis for our Specht
 modules  $S^{\Delta,\Delta^{'}}$, which correspond to the basis consisting
 of standard tableaux in the case of symmetric groups.

\noindent
{\bf Example 3.24} Let $\Phi={\bf D_{4}}$ 
 with simple system 
 \[\pi =\{\alpha_{1}=\epsilon_{1}-\epsilon_{2},
\alpha_{2}=\epsilon_{2}-\epsilon_{3},\alpha_{3}=\epsilon_{3}-\epsilon_{4},
\alpha_{4}=\epsilon_{3}+\epsilon_{4}\}\] 
\noindent
Let $~e~$, $~\tau_{2}~$,
 $(~\tau_{1}~\tau_{2}~)$,
 $(~\tau_{1}~\tau_{3}~\tau_{2}~)$, 
 $(~\tau_{1}~
\tau_{2}~\tau_{1}~\tau_{3}~
\tau_{2}~\tau_{1}~)$,
 $(~\tau_{1}~\tau_{4}~\tau_{2}~)$,
 $(~\tau_{1}~\tau_{2}~
\tau_{1}~\tau_{4}~\tau_{2}~
\tau_{1}~)$, $(~\tau_{3}~
\tau_{4}~\tau_{2}~)$, $(~\tau_{2}~\tau_{3}~\tau_{2}~
\tau_{4}~\tau_{2}~\tau_{3}~)$, $(~\tau_{1}~\tau_{3}~
\tau_{4}~\tau_{2}~)$, $(~\tau_{1}~\tau_{2}~\tau_{1}~
\tau_{3}~\tau_{2}~\tau_{1}~\tau_{4}~\tau_{2}~
\tau_{1}~\tau_{3}~\tau_{2}~\tau_{4}~)$,
$(~\tau_{2}~
\tau_{1}~\tau_{3}~\tau_{4}~\tau_{2}~)$,  $(~\tau_{2}~
\tau_{1}~\tau_{4}~\tau_{2}~
\tau_{1}~\tau_{3}~)$ be representatives
 of conjugate classes $C_{1}$,
 $C_{2}$, $C_{3}$, $C_{4}$, $C_{5}$, $C_{6}$, 
$C_{7}$, $C_{8}$, $C_{9}$, $C_{10}$,
 $C_{11}$, $C_{12}$, $C_{13}$ respectively of
 ${\cal W}({\bf D_{4}})$.  Then
the character table of  ${\cal W}({\bf D_{4}})$ is :

\vspace{0.2cm}
\begin{tabular}{r|rrrrrrrrrrrrr}
&$C_{1}$&$C_{2}$&$C_{3}$&$C_{4}$&$C_{5}$&$
C_{6}$&$C_{7}$&$C_{8}$&$C_{9}$&$C_{10}$&$
C_{11}$&$C_{12}$&$C_{13}$\\
\hline  
$\chi_{1}$&1&1&1&1&1&1&1&1&1&1&1&1&1\\ 
$\chi_{2}$&1&-1&1&-1&1&-1&1&-1&1&1&1&-1&1\\
$\chi_{3}$&2&0&-1&0&2&0&2&0&2&-1&2&0&2\\
$\chi_{4}$&3&1&0&1&3&-1&-1&-1&-1&0&3&1&-1\\
$\chi_{5}$&3&1&0&-1&-1&-1&-1&1&3&0&3&1&-1\\
$\chi_{6}$&3&-1&0&-1&3&1&-1&1&-1&0&3&-1&-1\\
$\chi_{7}$&3&-1&0&1&-1&1&-1&-1&3&0&3&-1&-1\\
$\chi_{8}$&3&-1&0&1&-1&-1&3&1&-1&0&3&-1&-1\\
$\chi_{9}$&3&1&0&-1&-1&1&3&-1&-1&0&3&1&-1\\
$\chi_{10}$&4&-2&1&0&0&0&0&0&0&-1&-4&2&0\\
$\chi_{11}$&4&2&1&0&0&0&0&0&0&-1&-4&-2&0\\
$\chi_{12}$&6&0&0&0&-2&0&-2&0&-2&0&6&0&2\\
$\chi_{13}$&8&0&-1&0&0&0&0&0&0&1&-8&0&0
\end{tabular}  

The non-conjugate subsystems of ${\bf D_{4}}$ are:

\begin{tabular}{lll}
(1)~~$\Psi_{1}~~=~~\bf{A_{3}}$\normalsize&  
with simple system ~$J_{1}~=$&$~\{1000,0100,0010\}~$
\\
(2)~~$\Psi_{2}~~=~~\bf{A_{3}^{'}}$\normalsize & 
with simple system~ $J_{2}~=$&$~\{1000,0100,0001\}~$
\\
(3)~~$\Psi_{3}~~=~~\bf{A_{3}^{''}}$\normalsize& 
 with simple system ~$J_{3}~=$&$~\{0100,0010,0001\}$
\\
(4)~~$\Psi_{4}~~=~~\bf{4A_{1}}$\normalsize  & 
with simple system ~$J_{4}~=$&$~\{1000,1211,0010,0001\}$
\\
(5)~~$\Psi_{5}~~=~~\bf{3A_{1}}$\normalsize  & 
with simple system ~$J_{5}~=$&$~\{1000,0010,0001\}$
\\
(6)~~$\Psi_{6}~~=~~\bf{A_{2}}$\normalsize & with 
simple system~ $J_{6}~=$&$~\{1000,0100\}$
\\
(7)~~$\Psi_{7}~~=~~\bf{2A_{1}}$\normalsize& 
 with simple system ~$J_{7}~=$&$~\{1000,0010\}$
\\
(8)~~$\Psi_{8}~~=~~\bf{2A_{1}^{'}}$\normalsize &
 with simple system ~$J_{8}~=$&$~\{1000,0001\}$
\\
(9)~~$\Psi_{9}~~=~~\bf{2A_{1}^{''}}$\normalsize & 
with simple system ~$J_{9}~=$&$~\{0010,0001\}$
\\
(10) $\Psi_{10}~=~\bf{A_{1}}$\normalsize&  with
 simple system $J_{10}~=$&$~\{1000\}$
\\
(11) $\Psi_{11}~=~\bf{D_{4}}$\normalsize&  with 
simple system $J_{11}~=$&$~\{1000,0100,0010,0001\}~$
\\
(12) $\Psi_{12}~=~\emptyset$&  with simple 
system $J_{12}~=$&$~\emptyset$ .
\end{tabular}

\vspace{0.1cm}

Let $\Psi_{1}={\bf A_{3}}$ be the subsystem 
of ${\bf D_{4}}$ with $J_{1}=\{1000,0100,0010\}$. Let $\Psi_{1}^{'}={\bf 2A_{1}}$ 
 be the subsystem of $\Phi$ which is 
contained in $\Phi \setminus \Psi$, with simple system $J_{1}^{'}=\{1101,0111\}$. Since ${\cal W}(J_{1}) \cap {\cal W}(J_{1}^{'})=<e>$ and   
${\cal W}(J_{1}^{\perp}) \cap {\cal W}(J_{1}^{'^{\perp}})=<e>$  , then $\{J_{1},J_{1}^{'}\}$ is a useful system in $\Phi$ . Then $\tau_{\Delta_{1}}$ contains the $\Delta_{1}-tabloids$ ;

\begin {tabular}{lll}
$\{ \bar{J_{1}} \}$&=&$\{1000,0100,0010;1101,0111\}~$\\
$\{ \tau_{4}\bar{J_{1}} \}$&=&$\{1000,0101,0010;1100,0110 \}~$\\
$\{ \tau_{2}\tau_{4}\bar{J_{1} }\}$&=&$\{1100,0001,0110;1000,0010\}~$\\
$\{ \tau_{1}\tau_{2}\tau_{4}\bar{J_{1}} \}$&=&$\{0100,0001,1110;-1000,0010\}~$\\
$\{ \tau_{3}\tau_{2}\tau_{4}\bar{J_{1}} \}$&=&$\{1110,0001,0100;1000,-0010\}~$\\
$\{ \tau_{1}\tau_{3}\tau_{2}\tau_{4}\bar{J_{1}} \}$&=&$\{0110,0001,1100;-1000,-0010\}~$\\
$\{ \tau_{2}\tau_{1}\tau_{3}\tau_{2}\tau_{4}\bar{J_{1}} \}$&=&$\{0010,0101,1000;-1100,-0110\}~$\\
$\{ \tau_{4}\tau_{2} \tau_{1}\tau_{3}\tau_{2}\tau_{4}\bar{J_{1}} \}$&=&$\{0010,0100,1000;-1101,-0111\}$
\end{tabular}

\vspace{0.1cm}
\noindent
For $d=e,\tau_{1}\tau_{2}\tau_{4},\tau_{3}\tau_{2}
\tau_{4},\tau_{4}\tau_{2}\tau_{1}\tau_{3}\tau_{2}
\tau_{4}$ we have  $d\Psi_{1}\cap\Psi_{1}^{'}=\emptyset$. Since
\[e_{J_{1},J_{1}^{'}}=\{\bar{J} \}-\{ \tau_{1}\tau_{2}
\tau_{4}\bar{J} \}-\{ \tau_{3}\tau_{2}\tau_{4}\bar{J} \}+\{ \tau_{4}\tau_{2} 
\tau_{1}\tau_{3}\tau_{2}\tau_{4}\bar{J} \}\]
\noindent
then $\{J_{1},J_{1}^{'}\}$ is a good system in $\Phi$ .

Now  let $K$ be a field with Char$K$=0. Let $M^{\Delta_{1}}$
 be $K$-space whose basis elements are the $\Delta_{1}$-tabloids. 
 Let $S^{\Delta_{1},\Delta_{1}^{'}}$ be the corresponding $K{\cal W}$-submodule of 
 $M^{\Delta_{1}}$, then  by definition of the Specht module we have 
\[S^{\Delta_{1},\Delta_{1}^{'}}=Sp~\{~e_{J_{1},J_{1}^{'}}~,~e_{\tau_{4}J_{1},\tau_{4}J_{1}^{'}}~,~e_{\tau_{2}\tau_{4}J_{1},
\tau_{2}\tau_{4}J_{1}^{'}}~\}\]

\noindent
where

$\begin{array}{lll}
e_{J_{1},J_{1}^{'}}&=&\{\bar{J_{1}} \}-\{ \tau_{1}\tau_{2}
\tau_{4}\bar{J}_{1} \}-\{ \tau_{3}\tau_{2}\tau_{4}\bar{J}_{1} \}+
\{ \tau_{4}\tau_{2} \tau_{1}\tau_{3}\tau_{2}\tau_{4}\bar{J}_{1} \}\\

e_{\tau_{4}J_{1},\tau_{4}J_{1}^{'}}&
=&\{\tau_{4}\bar{J}_{1} \}-\{ \tau_{1}\tau_{2}
\tau_{4}\bar{J} _{1}\}-\{ \tau_{3}\tau_{2}\tau_{4}\bar{J}_{1} \}+\{ \tau_{2} 
\tau_{1}\tau_{3}\tau_{2}\tau_{4}\bar{J}_{1} \}\\

e_{\tau_{2}\tau_{4}J_{1},\tau_{2}\tau_{4}J_{1}^{'}}&=&\{\tau_{2}\tau_{4}
\bar{J}_{1} \}-\{ \tau_{1}\tau_{2}
\tau_{4}\bar{J}_{1} \}-\{ \tau_{3}\tau_{2}\tau_{4}\bar{J}_{1} \}+
\{ \tau_{1}\tau_{3}\tau_{2}\tau_{4}\bar{J}_{1}\}
\end{array}$

Let $T_{1}$ be the matrix representation of ${\cal W}$ afforded by 
$S^{\Delta_{1},\Delta_{1}^{'}}$ with character $\psi_{1}$  and let $\tau_{2}$ be the representative of the
conjugate class $C_{2}$. Then

{\setlength{\arraycolsep}{0.1cm}
$\begin{array}{lllll}
\tau_{2}(e_{J_{1},J_{1}^{'}})&=&\{\bar{J_{1}} \}-\{ \tau_{1}\tau_{2}
\tau_{4}\bar{J}_{1} \}-\{ \tau_{3}\tau_{2}\tau_{4}\bar{J}_{1} \}+
\{ \tau_{4}\tau_{2} \tau_{1}\tau_{3}\tau_{2}\tau_{4}\bar{J}_{1} \}&=&
e_{J_{1},J_{1}^{'}}\\

\tau_{2}(e_{\tau_{4}J_{1},\tau_{4}J_{1}^{'}})&
=&\{\tau_{2}\tau_{4}\bar{J}_{1} \}-\{ \tau_{1}\tau_{2}
\tau_{4}\bar{J} _{1}\}-\{ \tau_{3}\tau_{2}\tau_{4}\bar{J}_{1} \}+
\{\tau_{1}\tau_{3}\tau_{2}\tau_{4}\bar{J}_{1} \}&=&
e_{\tau_{2}\tau_{4}J_{1},\tau_{2}\tau_{4}J_{1}^{'}}\\

\tau_{2}(e_{\tau_{2}\tau_{4}J_{1},\tau_{2}\tau_{4}J_{1}^{'}})&=&
\{\tau_{4}\bar{J}_{1} \}-\{ \tau_{1}\tau_{2}
\tau_{4}\bar{J} _{1}\}-\{ \tau_{3}\tau_{2}\tau_{4}\bar{J}_{1} \}+\{ \tau_{2} 
\tau_{1}\tau_{3}\tau_{2}\tau_{4}\bar{J}_{1} \}&=&e_{\tau_{4}J_{1},\tau_{4}J_{1}^{'}}
\end{array}$}

Thus we have

$T_{1}~(~\tau_{2}~)~=~\left( \begin{array}{ccc}
1&0&0\\
0&0&1\\
0&1&0
\end{array}
\right) $
     and $\psi_{1}(\tau_{2})=1$.

By a similar calculation to the above it can be showed that $\psi_{1}=\chi_{4}$. By the same method to the above, 
we have
\vspace{0.2cm}

{\setlength{\arraycolsep}{0.2cm}
$\begin{array}{l|l|l|l|l}
\Psi&~~~J&\Psi^{'}&~~~J^{'}&Ch\\
\hline
{\bf A_{3}}&\{1000,0100,0010\}&{\bf
2A_{1}}&\{1101,0111\}&\chi_{4}
\\
{\bf A_{3}^{'}}&\{1000,0100,0001\}&{\bf
2A_{1}^{'}}&\{1110,0111\}&\chi_{9}
\\
{\bf A_{3}^{''}}&\{0100,0010,0001\}&{\bf
2A_{1}^{''}}&\{1110,1101\}&\chi_{5}
\\
{\bf 4A_{1}}&\{1000,1211,0010,0001\}&
{\bf 3A_{1}}&\{1100,0110,0101\}&\chi_{3}
\\
{\bf 3A_{1}}&\{1000,0010,0001\}&
{\bf A_{2}}&\{1111,0100\}&\chi_{13}
\\
{\bf A_{2}}&\{1000,0100\}&{\bf 3A_{1}}&\{1110,1101,0111\}&
\chi_{13}
\\
{\bf 2A_{1}}&\{1000,0010\}&{\bf
A_{3}}&\{0001,1100,0110\}&\chi_{7}
\\
{\bf 2A_{1}^{'}}&\{1000,0001\}&{\bf
A_{3}^{'}}&\{0010,0101,1100\}&\chi_{6}
\\
{\bf 2A_{1}^{''}}&\{0010,0001\}&{\bf
A_{3}^{''}}&\{1000,0101,0110\}&\chi_{8}
\\
{\bf D_{4}}&\{1000,0100,0010,0001\}&
~~\emptyset&~~\emptyset&\chi_{1}
\\
~~\emptyset&~~\emptyset&
{\bf D_{4}}&\{1000,0100,0010,0001\}&\chi_{2}

\end{array}$}
\vspace{0.2cm}

We note that the irreducible modules corresponding  to the characters
$\chi_{10}$, $\chi_{11}$, $\chi_{12}$ have not been obtained. We now show how
an additional irreducible character is obtained. Let $\Psi_{2}^{'}={\bf
A_{1}}$ be the subsystem of $\Phi$ with simple system
$J_{2}^{'}=\{1101\}$. Then $\{J_{1},J_{2}^{'}\}$ is a useful system in $\Phi$.
Since $\Psi_{2}^{'} \subset  \Psi_{1}^{'}$, by Lemma 3.15 
$S^{\Delta_{1},\Delta_{1}^{'}}$ is a $K{\cal W}$-submodule of 
$S^{\Delta_{1},\Delta_{2}^{'}}$. By a similar calculation to the above it can
be showed that the corresponding character of ${\cal W}$ afforded by 
$S^{\Delta_{1},\Delta_{2}^{'}}/S^{\Delta_{1},\Delta_{1}^{'}}$ is $\chi_{11}$.

In fact, we have a '$Specht~series$' corresponding to series
\[ \emptyset \subset \{1101\} \subset \{1101,0111\}\]
\noindent
Unfortunately, no further irreducible character of 
${\cal W}({\bf D_{4}})$ are obtained by similar
calculations for the other subsystems.$~~~~~\Box$

We see that not all the irreducible
modules for ${\cal W}({\bf D_{4}})$ are obtained in this way. That is, 
 there are not a sufficient supply of subsystems or Weyl subgroups to give
 a complete set of irreducible modules. In the next section, it
will be shown that the irreducible character $\chi_{10},\chi_{12}$ of degree
4 and 6 respectively can be obtained by considering ${\cal W}({\bf
B_{3}})$ and  ${\cal W}({\bf G_{2}})$ as Steinberg subgroups of 
${\cal W}({\bf D_{4}})$.

\section{Additional Specht Modules for Weyl Groups}

In this section we consider our groups as reflection groups and thus,
 include reflection subgroups these groups. We modify the algorithm 
 for determining subsystems of $\Phi$ in ( 2.4 )  so as to include 
 Steinberg subgroups of the type 
described in ( 2.8 ). Detailed proofs are not always included as they are
either modifications of those in Section 3 or in the earlier paper [ 11 ].
In this case, it is necessary to express our results in terms of the 
subgroups rather than the root subsystems as in Section 3.

\noindent
{\bf Algorithm 4.1}
\newline
{\it (i) Form the extended Dynkin diagram of $\Phi$ . This is obtained by adding one
further node to the graph of $\Phi$ corresponding to the negative of the highest
root ,
\newline
(ii) Delete one or more nodes in all possible ways from the extended Dynkin 
diagram,
\newline
(iii) Let $\rho$ be a symmetry of the remaining Dynkin diagram.
Then follow (2.8)
to form a Steinberg subgroup for this remaining Dynkin diagram ,
\newline
(iv) Take also the duals of the diagrams obtained in the same way from the 
dual system $\tilde{\Phi}$ which is obtained from by interchancing long and 
short roots,
\newline
(v) Repeat the process with the Dynkin diagram obtained and continue any 
number of times.}$~\Box$

\noindent
{\bf Definition 4.2}  An $extended~ subsystem$ $\Psi$  is 
the root system  corresponding to the Dynkin diagram 
obtained by  Algorithm 4.1. $~~\Box$

Let $\Psi$  be an extended subsystem of $\Phi$ with simple
system $J$ and Dynkin diagram $\Delta$. By  Algorithm 4.1 we have
\newline
(i)  If the symmetry $\rho$ is the identity, then 
${\cal W}(J)$ is a Weyl subgroup
of ${\cal W}(\Phi)$. 
\newline
(ii) If $\rho$ is a non-trivial symmetry, then ${\cal W}(J)$ is a 
Steinberg subgroup of ${\cal W}(\Phi)$.

\noindent
{\bf Definition 4.3}
If $\Psi$ is an extended subsystem of $\Phi$, then the reflection 
subgroup ${\cal W}(J)$ of ${\cal W}$ which is generated by 
the $\tau_{\alpha}~,~\alpha \in \Psi$ is called a $Steinberg-Weyl~(S{\cal W})-subgroup$ of ${\cal W}$.$~~~\Box$

We now define equivalence relation on the elements of ${\cal W}$.
\[w^{'} \sim w ~~~~if~and~only~if~~~~w^{'}w^{-1} \in {\cal W}(J)~
for~w ~,~ w^{'} \in {\cal W}(\Phi) \]

\noindent
{\bf Definition 4.4}
 Let $w \in {\cal W}$ . The equivalence class 

\[ \{~w~{\cal W}(J)~\}~=~\{~w^{'} \mid w^{'}~\sim~w \} \] 

\noindent
is called a $\Delta-tabloid$. $~~~~~~~~~\Box$

Let $\tau_{\Delta}$ be the set of all $\Delta$-tabloids. In this 
new setting it is clear that the elements of $\tau_{\Delta}$ are the left coset 
  representatives of ${\cal W}(J)$ in ${\cal W}$ and the number of distinct 
 elements in $\tau_{\Delta}$ is $[{\cal W}~:{\cal W}(J)]$.

The group ${\cal W}$ acts on $\tau_{\Delta}$ according to 
\[w^{'}~\{~w~{\cal W}(J)~\} =\{w^{'}~w~{\cal W}(J)~\}~~~~~~for~all~w^{'} 
\in {\cal W} .\]
\noindent

This action is again easily seen to be well defined .

Now if K is arbitrary field , let $M^{\Delta}$ be the K-space whose basis elements are $\Delta$ -tabloids. Extend the action of ${\cal W}$ on $\tau_{\Delta}$ linearly on $M^{\Delta}$, then $M^{\Delta}$ becomes $K{\cal W}$-module.  
Then we have the following theorem.

\noindent
{\bf Theorem 4.5}
{\it The $K{\cal W}$-module  $M^{\Delta}$ is the permutation module on the
 $S{\cal W}$-subgroup ${\cal W}(J)$  $ M^{\Delta}$  is a cyclic 
$K{\cal W}$-module generated by any one tabloid and $dim_{K}M^{\Delta} =
 [~{\cal W}~:~{\cal W}(J)~]$.}

Now we can consider the possibility of constructing a $K{\cal W}$-module 
$S^{J,J^{'}}$ which generalises the Specht modules described in Section 3.
In this direction we first define a useful dual of a  $S{\cal W}$-subgroup ${\cal W}(J)$.

\noindent
{\bf Definition 4.6}
A $useful~ dual$ of ${\cal W}(J)$ is a $S{\cal W}$-subgroup  ${\cal W}(J^{'})$ of ${\cal W}(\Phi)$ which satisfies ${\cal W}(J)\cap~{\cal W}(J^{'}) = < ~e~>~~~~~~~~\Box$

\noindent
Then we have the following lemma.

\noindent
{\bf Lemma 4.7}
{\it  If  ${\cal W}(J^{'})$ is a useful dual of ${\cal W}(J)$ then 
${\cal W}(w~J^{'})$ is also a useful dual of ${\cal W}(J)$ for all $w \in {\cal W}(J)$.}

\noindent
{\bf Definition 4.8}
Let ${\cal W}(J)$ be a $S{\cal W}$-subgroup of ${\cal W}(\Phi)$  and  
${\cal W}(J^{'})$ a useful dual of ${\cal W}(J)$. Let
\begin {eqnarray*}
\kappa_{J^{'}}=\sum_{\sigma~\in~{\cal W}(J^{'})}~s~(~\sigma~)~\sigma~~~and~~~
e_{J,J^{'}}=~\kappa_{J^{'}}\{~{\cal W}(J)~\}
\end{eqnarray*}
\noindent
where $s$ is the sign function defined in Section 2.Then $e_{J,J^{'}}$ 
is called a $generalized$ $\Delta$ $-polytabloid$.

 Let $S^{J,J^{'}}$ be the subspace of $M^{\Delta}$ 
spanned by all the  generalized 
$\Delta$-polytabloids 
$e_{wJ,wJ^{'}}$,  $w\in~{\cal W}$. By the same method as in Section 3, 
$S^{J,J^{'}}$ is a 
$K{\cal W}$-submodule of $M^{\Delta}$ , 
which is called a $generalized$ $Specht~ module$.

If  ${\cal W}(J^{'})~=~<~e~>$  
  then  $S^{J,J^{'}} \cong M^{\Delta}$.

\noindent
{\bf Theorem 4.9} 
{\it  $S^{J,J^{'}}$ is a cyclic submodule generated by any $\Delta$-polytabloid.
}

\noindent
{\bf Lemma 4.10} 
{\it  Let ${\cal W}(J)$ be a $S{\cal W}$-subgroup  of ${\cal W}$ and 
let ${\cal W}(J^{'})$ be a useful dual of ${\cal W}(J)$. Let  $w^{'}$ be a left coset representative of ${\cal W}(J^{'})$ in {\cal W}.
Then $S^{J,J^{'}}$ is spanned by $w^{'}~e_{J,J^{'}}$.}

\noindent
{\bf Lemma 4.11} 
{\it  Let ${\cal W}(J)$ be a $S{\cal W}$-subgroup of ${\cal W}$.  
If $\{w^{'}~{\cal W}(J) \}$ appears in $e_{J,J^{'}}$ then it appears only once.}

\noindent
{\bf Proof} See Lemma 3.11.$~~\Box$

\noindent
{\bf Corollary 4.12} 
{\it  If ${\cal W}(J) \cap {\cal W}(J^{'})=~<~e~>$ , then $e_{J,J^{'}} \neq 0$}.

\noindent
{\bf Lemma 4.13} 
{\it Let ${\cal W}(J_{1}^{'})$ and ${\cal W}(J_{2}^{'})$ be useful duals of 
${\cal W}(J)$. If ${\cal W}(J_{1}^{'} )$ is a $S{\cal W}$-subgroup of ${\cal W}( J_{2}^{'})$, then $S^{J,J_{2}^{'}}$ is a $K{\cal W}$-submodule of $S^{J,J_{1}^{'}}$.}

\noindent
{\bf Proof} 
It follows from Lemma 3.15 $\Box$

\noindent
{\bf Lemma 4.14} 
{\it If $w \in {\cal W}$ and ${\cal W}(J^{'})$ is a useful dual of ${\cal W}(J)$,  then the following conditions are equivalent:
\newline
(i) $\{~w~{\cal W}(J)~\}$ appears with non-zero coefficient in $e_{J,J^{'}}$
\newline
(ii) There exists $\sigma \in {\cal W}(J^{'})$ such that $\sigma \{~{\cal W}(J)~\}~=~\{~w~{\cal W}(J)~\}$
\newline
(iii) There exist $\rho \in {\cal W}(J)$, $\sigma \in {\cal W}(J^{'})$ such that $w~=~\sigma~\rho$.}

\noindent
{\bf Proof} 
 The equivalence of (i) and (ii) follows directly from the formula 
\[e_{J,J^{'}}=\sum_{\sigma~\in~{\cal W}(J^{'})}~s~(~\sigma~)~\sigma~\{~{\cal W}(J)~\}\]

\noindent
(ii) $\Rightarrow$ (iii) Suppose that there exists $\sigma
 \in {\cal W}(J^{'})$ such that  $\sigma \{{\cal W}(J)\}=\{w{\cal W}(J)\}$.
 Then we have $\sigma^{-1} w \{{\cal W}(J)\}=\{{\cal W}(J)\}$. 
 By the definition of 
equivalence,  $\sigma^{-1}w \in {\cal W}(J)$ and there exists $\rho \in
 {\cal W}(J)$ such that $\sigma^{-1} w=\rho$.  
 Hence $w=\sigma\rho$ , where $\sigma \in {\cal W}(J^{'})$ and
$\rho \in {\cal W}(J)$.
\newline
(iii) $\Rightarrow$ (ii) If $w=\sigma\rho$, then since $\rho \in {\cal W}(J)$, 
$\rho\{{\cal W}(J)\}=\{{\cal W}(J)\}$ and $ \{w{\cal W}(J)\}=\{\sigma{\cal W}(J)
\}.~\Box $

\noindent
{\bf Definition 4.15} 
 A useful dual ${\cal W}(J^{'})$  of  ${\cal W}(J)$  is called a $good~dual$ 
 of ${\cal W}(J)$ if $\kappa_{J^{'}}\{w{\cal W}(J)\}\neq 0$, then $\{w{\cal W}(J)\}$ appears in $e_{J,J^{'}}$.$~~~~~~~~\Box$

\noindent
{\bf Lemma 4.16} 
{\it  Let $~{\cal W}(J^{'})$ be a good dual of ${\cal W}(J)$ .
\newline
(i) If $\{w{\cal W}(J)\}$ does not appear in $e_{J,J^{'}}$  
then $\kappa_{J^{'}}\{w{\cal W}(J)\}=0.$
\newline
(ii) If $\{w{\cal W}(J)\}$ appears in $e_{J,J^{'}}$  then there exists $\sigma \in {\cal W}(J^{'})$ such that}
\[\kappa_{J^{'}}\{~w~{\cal W}(J)~\}~=~s~(~\sigma~)~e_{J,J^{'}} \]

\noindent
{\bf Proof} 
See Lemma 3.20 $\Box$

\noindent
{\bf Corollary 4.17}{\it If $m \in M^{\Delta}$ then $\kappa_{J^{'}}~m$ is a multiple of $e_{J,J^{'}}$.}

We now define a bilinear form $ <~,~>$ on $M^{\Delta}$ by setting
\[<\{w_{1}{\cal W}(J)\},\{w_{2}{\cal W}(J)\}>=\left\{
\begin{array}{ll} 1 &\mbox{$if~w_{1}~=~w_{2}$}\\
0&\mbox{otherwise} 
\end{array}
\right.
\]

\noindent
This is a symmetric, non-singular, ${\cal W}$-invariant, bilinear form on $M^{\Delta}$.

Now the analogue of James' submodule theorem can be proved in this more 
general setting.

\noindent
{\bf Theorem 4.18} 
{\it  Let U be submodule of $M^{\Delta}$ . Then either 
\newline
$S^{J,J^{'}} \subseteq U$ or $U \subseteq S^{J,J^{'^{\perp}}}$ where $S^{J,J^{'^{\perp}}}$ is complement of $S^{J,J^{'}}$ in $M^{\Delta}$ }.

We can now prove our principal result.

\noindent
{\bf Theorem 4.19} 
{\it The K{\cal W}-module $D^{J}~=~S^{J,J^{'}}~/~S^{J,J^{'}}~\cap S^{J,J^{'^{\perp}}}~$ is zero or irreducible.}

\noindent
{\bf Proof} 
See Theorem 3.23 $\Box$

Now to illustrate the above we show how a complete system of irreducible 
modules is determined in the case $\Phi={\bf D_{4}}$.

\noindent
{\bf Example 4.20} Let $\Phi = {\bf D_{4}}$ with 
$\pi =\{\alpha_{1}=\epsilon_{1}-\epsilon_{2}, \alpha_{2}=
\epsilon_{2}-\epsilon_{3},\alpha_{3}=\epsilon_{3}-\epsilon_{4},
\alpha_{4}=\epsilon_{3}+ \epsilon_{4}~\}$ and let $\Psi_{1}={\bf A_{3}}$ be the extended subsystem 
of ${\bf D_{4}}$ with  simple system $J_{1}=\{1000,0100,0010\}$.
Let $\Psi_{1}^{'}={\bf 2A_{1}}$ 
 be the extended subsystem of ${\bf D_{4}}$ with simple system 
 $J_{1}^{'}=\{1101,0111\}$.Since 
 ${\cal W}(J_{1}) \cap {\cal W}(J_{1}^{'})=<e>$, then ${\cal W}(J_{1}^{'})$ is 
a useful dual of ${\cal W}(J_{1})$. Then $\tau_{\Delta_{1}}$ contains
 the $\Delta_{1}$-tabloids $\{ {\cal W}(J_{1}) \}$,$\{ \tau_{4}{\cal W}(J_{1}) \}$,
$\{ \tau_{2}\tau_{4}
{\cal W}(J_{1})\}$, $\{ \tau_{1}\tau_{2}\tau_{4}{\cal W}(J_{1}) \}$, 
$\{ \tau_{3}\tau_{2}\tau_{4}{\cal W}(J_{1}) \}$, 
$\{ \tau_{1}\tau_{3}\tau_{2}\tau_{4}{\cal W}(J_{1}) \}$, 
$\{ \tau_{2}\tau_{1}\tau_{3}\tau_{2}\tau_{4}{\cal W}(J_{1}) \}$, 
$\{ \tau_{4}\tau_{2} \tau_{1}\tau_{3}\tau_{2}\tau_{4}{\cal W}(J_{1}) \}$.

For $w=e,\tau_{1}\tau_{2}\tau_{4},\tau_{3}\tau_{2}
\tau_{4},\tau_{4}\tau_{2}\tau_{1}\tau_{3}\tau_{2}
\tau_{4}$, we have  $\kappa_{J_{1}^{'}}\{w{\cal W}(J_{1})\}\neq 0$. Since
\[e_{J_{1},J_{1}^{'}}=\{{\cal W}(J_{1}) \}-\{ \tau_{1}\tau_{2}
\tau_{4}{\cal W}(J_{1}) \}-\{ \tau_{3}\tau_{2}\tau_{4}{\cal W}(J_{1}) \}+
\{ \tau_{4}\tau_{2} 
\tau_{1}\tau_{3}\tau_{2}\tau_{4}{\cal W}(J_{1}) \}\]
\noindent
then ${\cal W}(J_{1}^{'})$ is a good dual of ${\cal W}(J_{1})$.

Now  let $K$ be a field with Char$K$=0. Let $M^{\Delta_{1}}$
 be $K$-space whose basis elements are the $\Delta_{1}$-tabloids. 
 Let $S^{J_{1},J_{1}^{'}}$ be the corresponding $K{\cal W}$-submodule of 
 $M^{\Delta_{1}}$, then  by definition of the Specht module we have 
\[S^{J_{1},J_{1}^{'}}=Sp~\{~e_{J_{1},J_{1}^{'}}~,~e_{\tau_{4}J_{1},\tau_{4}J_{1}^{'}}~,~e_{\tau_{2}\tau_{4}J_{1},
\tau_{2}\tau_{4}J_{1}^{'}}~\}\]

\noindent
where

$\begin{array}{lll}
e_{J_{1},J_{1}^{'}}&=&\{{\cal W}(J_{1}) \}-\{ \tau_{1}\tau_{2}
\tau_{4}{\cal W}(J_{1}) \}-\{ \tau_{3}\tau_{2}\tau_{4}{\cal W}(J_{1}) \}+
\{ \tau_{4}\tau_{2} \tau_{1}\tau_{3}\tau_{2}\tau_{4}{\cal W}(J_{1}) \}\\

e_{\tau_{4}J_{1},\tau_{4}J_{1}^{'}}&
=&\{\tau_{4}{\cal W}(J_{1}) \}-\{ \tau_{1}\tau_{2}
\tau_{4}{\cal W}(J_{1})\}-\{ \tau_{3}\tau_{2}\tau_{4}{\cal W}(J_{1}) \}+
\{ \tau_{2} 
\tau_{1}\tau_{3}\tau_{2}\tau_{4}{\cal W}(J_{1}) \}\\

e_{\tau_{2}\tau_{4}J_{1},\tau_{2}\tau_{4}J_{1}^{'}}&=&\{\tau_{2}\tau_{4}
{\cal W}(J_{1}) \}-\{ \tau_{1}\tau_{2}
\tau_{4}{\cal W}(J_{1}) \}-\{ \tau_{3}\tau_{2}\tau_{4}{\cal W}(J_{1}) \}+
\{ \tau_{1}\tau_{3}\tau_{2}\tau_{4}{\cal W}(J_{1})\}
\end{array}$

Let $T_{1}$ be the matrix representation of ${\cal W}$ afforded by 
$S^{J_{1},J_{1}^{'}}$ with character $\psi_{1}$  and let $\tau_{2}$ be the representative of the
conjugate class $C_{2}$. Then

{\setlength{\arraycolsep}{0.01cm}
$\begin{array}{llll}
\tau_{2}(e_{J_{1},J_{1}^{'}})&=&\{{\cal W}(J_{1}) \}-\{ \tau_{1}\tau_{2}
\tau_{4}{\cal W}(J_{1}) \}-\{ \tau_{3}\tau_{2}\tau_{4}{\cal W}(J_{1}) \}+
\{ \tau_{4}\tau_{2} \tau_{1}\tau_{3}\tau_{2}\tau_{4}{\cal W}(J_{1}) \}&\\
&=&e_{J_{1},J_{1}^{'}}&\\

\tau_{2}(e_{\tau_{4}J_{1},\tau_{4}J_{1}^{'}})&
=&\{\tau_{2}\tau_{4}{\cal W}(J_{1}) \}-\{ \tau_{1}\tau_{2}
\tau_{4}{\cal W}(J_{1})\}-\{ \tau_{3}\tau_{2}\tau_{4}{\cal W}(J_{1}) \}+
\{\tau_{1}\tau_{3}\tau_{2}\tau_{4}{\cal W}(J_{1}) \}&\\
&=&e_{\tau_{2}\tau_{4}J_{1},\tau_{2}\tau_{4}J_{1}^{'}}&\\

\tau_{2}(e_{\tau_{2}\tau_{4}J_{1},\tau_{2}\tau_{4}J_{1}^{'}})&=&
\{\tau_{4}{\cal W}(J_{1}) \}-\{ \tau_{1}\tau_{2}
\tau_{4}{\cal W}(J_{1})\}-\{ \tau_{3}\tau_{2}\tau_{4}{\cal W}(J_{1}) \}+
\{ \tau_{2} 
\tau_{1}\tau_{3}\tau_{2}\tau_{4}{\cal W}(J_{1}) \}&\\
&=&e_{\tau_{4}J_{1},\tau_{4}J_{1}^{'}}&
\end{array}$}

Thus we have

$T_{1}~(~\tau_{2}~)~=~\left( \begin{array}{ccc}
1&0&0\\
0&0&1\\
0&1&0
\end{array}
\right) $
     and $\psi_{1}(\tau_{2})=1$.

By a similar calculation to the above it can be showed that 
$\psi_{1}=\chi_{4}$. By the same method to the above, 
in addition to the irreducible modules in Example 3.24 we have,
\vspace{0.2cm}

{\setlength{\arraycolsep}{0.3cm}
$\begin{array}{l|l|l|l|l}
\Psi&~~~J&\Psi^{'}&~~~J^{'}&Char\\
\hline
{\bf G_{2}}&\{0100,\frac{1}{3}(1000+0010+0001)\}&{\bf A_{2}}&\{0001,0110\}&
\chi_{12}
\end{array}$}
\vspace{0.2cm}

We note that the irreducible modules corresponding  to the characters
$\chi_{10},\chi_{11}$ have not been obtained. We now show how
 additional irreducible characters are obtained.

 Let $\Psi_{2}^{'}={\bf
A_{1}}$ be the extended subsystem of $\Phi$ with simple system
$J_{2}^{'}=\{1101\}$. Then ${\cal W}(J_{2}^{'})$ is a useful dual of 
${\cal W}(J_{1})$. Since ${\cal W}(J_{2}^{'})$ is a $S{\cal W}$-subgroup 
of ${\cal W}(J_{1}^{'})$, by Lemma 4.13 
$S^{J_{1},J_{1}^{'}}$ is a $K{\cal W}$-submodule of 
$S^{J_{1},J_{2}^{'}}$. By a similar calculation to the above it can
be showed that the corresponding character of ${\cal W}$ afforded by 
$S^{J_{1},J_{2}^{'}}/S^{J_{1},J_{1}^{'}}$ is $\chi_{11}$.

 Let $\Psi_{2}={\bf
A_{1}}$ be the extended subsystem of  $\Phi$ with simple system
$J_{2}=\{1000\}$. Let $\Psi_{1}^{'}={\bf B_{3}}$ be the extended subsystem 
of $\Phi$ with simple system
$J_{1}^{'}=\{0100,0001,\frac{1}{2}(1000+1211)\}$ and  $\Psi_{2}^{'}={\bf A_{3}}$ be another extended subsystem 
of $\Phi$ with simple system
$J_{2}^{'}=\{0100,0001,0010\}$. Then ${\cal W}(J_{2}^{'})$ is a useful dual 
of ${\cal W}(J_{2})$ and ${\cal W}(J_{1}^{'})$ is a good dual of 
${\cal W}(J_{1})$. Since ${\cal W}(J_{2}^{'})$ is a 
$S{\cal W}$-subgroup of ${\cal W}(J_{1}^{'})$, by Lemma 4.13 
$S^{J_{2},J_{1}^{'}}$ is a $K{\cal W}$-submodule of 
$S^{J_{2},J_{2}^{'}}$. By a similar calculation to the above it can
be showed that the corresponding character of ${\cal W}$ afforded by 
$S^{J_{2},J_{2}^{'}}/S^{J_{2},J_{1}^{'}}$ is $\chi_{10}$. 
Thus we have obtained a complete set of irreducible modules for ${\bf
D_{4}}$.$~\Box$
 
{\parbox{7cm}{\centering
Department of Mathematics\\
Ankara University\\
06100  Tando{\u g}an  Ankara\\
Turkey\\}
\parbox{7cm}{\centering
Department of Mathematics\\
The University of Wales\\
Aberystwyth SY23 3BZ\\
United Kingdom}
\end{document}